\documentstyle{article}
\def\zb{\overline{z}}
\def\Tb{\overline{T}}
\def\nn {\noindent}

\newtheorem{proposition}{Proposition}[section]

\newtheorem{lemma}[proposition]{Lemma}
\begin{document}

\title{UNITARY REPRESENTATIONS OF SU\boldmath$_q$(2) ON THE PLANE FOR $q
\in \mbox{R}^+$ OR GENERIC $q \in S^1$\footnote{Presented at the 7th
Colloquium``Quantum Groups and Integrable Systems'', Prague, 18--20 June 1998}}
\author{M. Irac-Astaud\footnote{E-mail: mici@ccr.jussieu.fr},\\
{\small \sl Laboratoire de Physique Th\'eorique de la Mati\`ere Condens\'ee,
Universit\'e Paris VII,}\\
{\small \sl 2,~place Jussieu, F-75251 Paris Cedex 05, France}\\
C. Quesne\footnote{Directeur de recherches FNRS;
E-mail: cquesne@ulb.ac.be},\\
{\small \sl Physique Nucl\'eaire Th\'eorique et Physique Math\'ematique,
Universit\'e Libre de Bruxelles,}\\ 
{\small \sl Campus de la Plaine CP229, Boulevard du
Triomphe, B-1050 Brussels, Belgium}}
\maketitle

\begin{abstract}
\nn Some time ago, Rideau and Winternitz introduced a realization of the quantum
algebra~su$_q$(2) on a real two-dimensional sphere, or a real plane, and
constructed a basis for its representations in terms of $q$-special
functions, which
can be expressed in terms of $q$-Vilenkin functions.
In their study, the values of~$q$ were implicitly restricted to $q \in
\mbox{R}^+$.
In the present paper, we extend their work to the case of generic values of
$q \in
S^1$ (i.e., $q$~values different from a root of unity). In
addition, we unitarize the representations for both types of $q$~values, $q \in
\mbox{R}^+$ and generic $q \in S^1$,  by determining some appropriate
scalar products.
\end{abstract}

%
%
\section{Introduction} \label{sec:intro}
 Rideau and Winternitz introduced a realization
of the
quantum algebra~su$_q$(2) on the plane~\cite{rideau} and constructed a
basis for its irreducible
representations (irreps) in terms of $q$-Vilenkin functions, so called because,
 for $q=1$,
they reduce to functions  introduced by
Vilenkin~\cite{vilenkin1}, and related to Jacobi polynomials.
%
%
This realization  was
used  to set up su$_q$(2)-invariant
Schr\"odinger equations in the usual framework of quantum
mechanics~\cite{irac1}.
%
%
Although not explicitly stated in Ref.~\cite{rideau}, the values of the
deformation
parameter~$q$, considered there, are restricted to $q \in \mbox{R}^+$.
Though important both from the $q$-special function viewpoint, and from that of
their applications in quantum mechanics, the question of the su$_q$(2) irrep
unitarity was also left unsolved.\par
%
%
The purpose of a previous paper \cite{nous} and of the present lecture
 is twofold: firstly, to find basis functions
 of the representation for generic $q \in S^1$ (i.e., for $q$ different
from a root
of unity),
and secondly, to unitarize the representations for both $q\in \mbox{R}^+$, and
generic $q
\in S^1$. As a consequence, the orthonormality relations of the
$q$-Vilenkin
and related functions are established~\cite{nous}.
%
In Sec.~\ref{sec:representations}, the representations of~su$_q$(2)
obtained in Ref.~\cite{rideau} are briefly reviewed. Their basis functions
are determined in Sec.~\ref{sec:Q}. The unitarization of the representations
is dealt
with in Sec.~\ref{sec:unitarization}. Sec.~\ref{sec:conclusion} contains the
conclusion.

%
\section{ Representations of su$_q$(2) on the plane}
\label{sec:representations}
The su$_q$(2) generators $H_3$, $H_+$, $H_-$ satisfy the commutation
relations~\cite{chari}
\begin{equation}
  \left[H_3, H_{\pm}\right] = \pm H_{\pm}, \qquad \left[H_+, H_-\right] =
  \left[2H_3\right]_q \equiv \frac{q^{2H_3} - q^{-2H_3}}{q - q^{-1}},
\label{eq:com}
\end{equation}
and the Hermiticity properties $ H_3^{\dagger} = H_3$, and $ H_{\pm}^{\dagger}
= H_{\mp}$.
One can
construct a Casimir operator $
  {\cal C} = H_+ H_- + \left[H_3\right]_q \left[H_3-1\right]_q = H_- H_+
  + \left[H_3\right]_q \left[H_3+1\right]_q$.
The generators $H_3$, $H_+$, $H_-$ can be realized by the
following
operators, acting on the plane, more precisely, on functions $f(z,\zb)$ of a
complex variable $z$ and its complex conjugate $\zb$ ~\cite{rideau},
\begin{eqnarray}
  H_3 & = & - T+ \Tb - N,
          \nonumber \\
  H_+ & = & - z^{-1}\, [T]_q\, q^{\Tb - (N/2)} - q^{T + (N/2)}\, \zb \left[\Tb -
          N\right]_q, \nonumber \\
  H_- & = & z\, [T + N]_q\, q^{\Tb - (N/2)} + q^{T + (N/2)}\, \zb^{-1}
          \left[\Tb\right]_q,  \label{eq:su-q}
\end{eqnarray}
where
$T = z \partial_z$, and $\Tb = \zb \partial_{\zb}$.\par
%
%
Basis functions $\Psi^J_{MNq}(z,\zb)$ for the ($2J+1$)-dimensional irrep
of~su$_q$(2) satisfy the relations~\cite{chari}
\begin{eqnarray}
  H_3 \Psi^J_{MNq} & = & M \Psi^J_{MNq}, \quad H_{\pm} \Psi^J_{MNq} = \left(
         [J\mp M]_q [J\pm M+1]_q\right)^{1/2} \Psi^J_{M\pm1,Nq}, \nonumber \\
  {\cal C} \Psi^J_{MNq} & = & [J]_q [J+1]_q \Psi^J_{MNq},\quad  M = \{ -J,-J+1,
         \cdots,J\} , \quad |N| \leq J,   \label{eq:irrep}
\end{eqnarray}
where $J$, $M$ and $N$ are simultaneously integers or half-integers. Let us
remark that, when $q \in S^1$, the existence of such a representation
implies that $[n]_q$ does not vanish unless $n=0$,
 hence that $q$ is not a root of unity.\par
%
%
{}Following Ref.~\cite{rideau}, let us write $ \Psi^J_{MNq}(z,\zb) $ as
\begin{equation}
 \Psi^J_{MNq}(z,\zb) = N^J_{MNq} Q_{Jq}(\eta)\, q^{-NM/2} R^J_{MNq}(\eta)\,
  \zb^{M+N}, \qquad \eta = z \zb.   \label{eq:Psi-z}
\end{equation}
Here, $N^J_{MNq}$ is a constant, which can be expressed as
\begin{equation}
  N^J_{MNq}  =  \frac{1}{\sqrt{2\pi}} \left(\frac{[J+N]_q!\,
          [2J+1]_q!}{[J-N]_q!}\right)^{1/2}    \left(\frac{[J+M]_q!}{[J-M]_q!\,
[2J]_q!}\right)^{1/2}\gamma(J,N,q),
\end{equation}
in terms of some yet undetermined normalization constant $\gamma(J,N,q)$. The
$q$-factorials are defined by $[x]_q! \equiv [x]_q [x-1]_q \ldots [1]_q$ if $x
\in \mbox{N}^+$,
$[0]_q! \equiv 1$, and $\left([x]_q!\right)^{-1} \equiv 0$ if $x \in
\mbox{N}^-$.
 The function $R^J_{MNq}(\eta)$ involved in Eq.~(\ref{eq:Psi-z})
is a polynomial
\begin{equation}
  R^J_{MNq}(\eta)  =
    \sum_k \frac{ [J-N]_q!\, [J-M]_q!(-\eta)^k}{[k]_q!\, [J-M-k]_q!\,
[J-N-k]_q!\,
          [M+N+k]_q!},   \label{eq:R}
\end{equation}
the summation over $k$ being restricted by the condition that all the
arguments of the
factorials in
the denominator be positive. The function  $Q_{Jq}(\eta)$  involved in
Eq.~(\ref{eq:Psi-z})
is defined by the functional equation
\begin{equation}
  Q_{Jq}(q^2\eta)(1+\eta) =  Q_{Jq}(\eta) (1+q^{-2J}\eta),      \label{eq:equQ}
\end{equation}
whose solution, only determined up to an arbitrary multiplicative factor
$f_{Jq}(\eta)$ such that
$f_{Jq}(q^2\eta) =  f_{Jq}(\eta)$,
will be discussed
in the
next section.\par
%
%
The functions~(\ref{eq:Psi-z}) are related to the $q$-Vilenkin functions,
defined by
\begin{equation}
  P^J_{MNq}(\xi)  =  i^{2J-M-N} \left(\frac{[J+M]_q!\,
[J+N]_q!}{[J-M]_q!\, [J-N]_q!}
          \right)^{1/2} \eta^{(M+N)/2} Q_{Jq}(\eta)\, R^J_{MNq}(\eta),
          \label{eq:q-Vil}
\end{equation}
where $ \eta = (1+\xi)/(1-\xi)$.
 For integer $J$~values, the functions
$\Psi^J_{M0q}$ are proportional to $q$-spherical harmonics, while
$P_{Jq}(\xi) \equiv P^J_{M0q}(\xi)$ are $q$-analogues of Legendre
polynomials.\par
%
%
In the $q \to 1$ limit, the su$_q$(2) realization~(\ref{eq:su-q}) goes
into a standard
su(2) realization, and choosing $Q_{J1}(\eta)  =
(1+\eta)^{-J}$, the functions $\Psi^J_{MN1}$ form an orthonormal set with
respect to the scalar product
\begin{equation}
  \langle \psi_1|\psi_2 \rangle = 2 \int \frac{dz d\zb}{(1+z\zb)^2}\,
  \overline{\psi_1\left(z,\zb\right)}\, \psi_2\left(z,\zb\right).
  \label{eq:prodscal}
\end{equation}

%
%
\section{Determination of $Q_{Jq}(\eta)$}
\label{sec:Q}
{}Following Ref.~\cite{rideau}, as a solution of
Eq.~(\ref{eq:equQ}),
we may consider the function
\begin{equation}
 Q_{Jq}(\eta)  =  {}_1\Phi_0\left(q^{2J}; -; q^2, -q^{-2J}\eta\right) =
          {}_1\Phi_0\left(q^{-2J}; -; q^{-2}, -q^{-2}\eta\right),
\label{eq:Q}
\end{equation}
where ${}_1\Phi_{0}$ is a basic hypergeometric series in the notations of
Ref.~\cite{exton}.\par
%
%
$\bullet$ For $q \in \mbox{R}^+$, use of the $q$-binomial theorem~\cite{exton}
leads to the expressions
\begin{equation}
\begin{array}{llll}
  Q_{Jq}(\eta)& =& \prod_{k=0}^{\infty}
\frac{(1+q^{2k}\eta)}{(1+q^{-2J+2k}\eta)},
 & \mbox {if $ 0 < q < 1$},
  \\
&&\\
 Q_{Jq}(\eta)& =& \prod_{k=0}^{\infty} \frac{(1+q^{-2J-2k-2}\eta)}
  {(1+q^{-2k-2}\eta)}, &\mbox{if $q > 1$}.
\end{array}
 \label{eq:QQ}
\end{equation}
For integer $J$~values, both expressions reduce to the inverse of a
polynomial,
\begin{equation}
  Q_{Jq}(\eta) = \prod_{k=0}^{J-1}\left( \frac{1}{1+\eta q^{-2J+2k}}\right).
  \label{eq:QJentier}
\end{equation}

\nn  For half-integer $J$~values, we are left with convergent infinite
products.\par
%
%
$\bullet$ For generic $q \in S^1$,
and integer $J$~values, Eq.~(\ref{eq:QJentier}) still
remains a valid solution of Eq.~(\ref{eq:equQ}). For half-integer
$J$~values, however, the infinite products contained in Eq.~(\ref{eq:QQ})
 are
divergent. We
have therefore to look for another solution to Eq.~(\ref{eq:equQ}).
 For such a purpose, by setting $
  K_{Jq}(\eta) = \ln Q_{Jq}(\eta)$,
we linearize Eq.~(\ref{eq:equQ}) into
\begin{equation}
  K_{Jq}(q^2 \eta) - K_{Jq}(\eta) = \ln \frac{1 + q^{-2J}\eta}{1 + \eta},
  \label{eq:equK}
\end{equation}
%
%
whose solution can be written in the form
$K_{Jq}(\eta) = L_q\left(q^{-2J-1}\eta\right) - L_q(q^{-1}\eta)$,
 where $L_q$ is solution of
\begin{equation}
 L_q(q\eta) - L_q(q^{-1}\eta) =
  \ln(1 + \eta).
  \label{eq:equL}
\end{equation}
In Ref.~\cite{nous}, we demonstrated
%
\begin{lemma}    \label{lem-repint}
{}For $0 < \eta < \infty$, and $q = e^{i\tau}$ different from a root of
unity, the
function
\begin{eqnarray}
  L_q(\eta) & = & \frac{1}{2\pi i} \int_0^{\infty} \frac{dt}{t(1+t)} \ln
\left(1 + \eta
          t^{\tau/\pi}\right), \qquad \mbox{if $0 < \tau < \pi$},
\label{eq:Lsol1} \\
  L_q(\eta) & = & - \frac{1}{2\pi i} \int_0^{\infty} \frac{dt}{t(1+t)} \ln
\left(1 +
          \eta t^{-\tau/\pi}\right), \qquad \mbox{if $-\pi < \tau < 0$},
\label{eq:Lsol2}
\end{eqnarray}
is a solution of Eq.~(\ref{eq:equL}).
\end{lemma}

%
%
The results of the present section can be collected into
%
\begin{proposition}
The function $Q_{Jq}(\eta)$, appearing on the right-hand side of
Eq.~(\ref{eq:Psi-z}),
is given by Eq.~(\ref{eq:QJentier}) for integer $J$~values, and either $q \in
\mbox{R}^+$ or generic $q \in S^1$, and by Eq.~(\ref{eq:QQ})  for
half-integer
$J$~values and $q \in \mbox{R}^+$. For half-integer $J$~values, and generic $q
\in S^1$,
it can be expressed as
\begin{equation}
  Q_{Jq}(\eta) = \exp\left\{L_q\left(q^{-2J-1}\eta\right) -
  L_q\left(q^{-1}\eta\right)\right\},  \label{eq:QJdemi-ent}
\end{equation}
where $L_q(\eta)$ admits the integral representation given in
Lemma~\ref{lem-repint}.
\end{proposition}

%
%
\section{Unitarization of the representations }
\label{sec:unitarization}
In the present section, we will determine a new scalar product $\langle \psi_1|
\psi_2 \rangle_q$ that unitarizes the realization~(\ref{eq:su-q}) of~su$_q$(2),
and goes over into the old one $\langle \psi_1| \psi_2 \rangle$, defined in
Eq.~(\ref{eq:prodscal}), whenever $q \to 1$.

%
%
\subsection{The case where \boldmath $q \in \mbox{R}^+$}
Let us make the following ansatz for $\langle \psi_1| \psi_2 \rangle_q$,
\begin{eqnarray}
  \langle \psi_1| \psi_2 \rangle_q & = & \int_0^{\infty} d\rho
\int_0^{2\pi} d\phi
          \biggl( \overline{A_q \psi_1(\rho,\phi,q)}\, f_1(\rho,q)\, q^{a_1 \rho
          \partial_{\rho}} \psi_2(\rho,\phi,q) \nonumber \\
  & & \mbox{} + \overline{\psi_1(\rho,\phi,q)}\, f_2(\rho,q)\, q^{a_2 \rho
          \partial_{\rho}} A_q \psi_2(\rho,\phi,q)\biggr), \qquad z=\rho
e^{i\phi}.
\label{eq:ansatz}
\end{eqnarray}
 Here $a_1$,~$a_2$ and $f_1(\rho,q)$, $f_2(\rho,q)$
are some
yet undetermined constants and functions, and $A_q \equiv q^{-2q \partial_q}$ is
the operator that changes $q$
into~$q^{-1}$, when acting on any function of~$q$.\par
%
%
It is easy to check that $ H_3$ satisfies the equation
 $\langle \psi_1| H_3 \psi_2 \rangle_q = \langle H_3 \psi_1| \psi_2 \rangle_q $.
Let us impose the condition
  $\langle \psi_1| H_+ \psi_2 \rangle_q = \langle H_- \psi_1| \psi_2 \rangle_q.$
We obtain~\cite{nous}
\begin{equation}
  a_1 = - 1, \quad a_2 = 1,
\label{eq:a}
\end{equation}
and
\begin{equation}
 f_1(\rho,q) = \frac{B_1(q) q^{-1} \rho}{\left(1+\rho^2\right)
\left(1+q^{-2}\rho^2
  \right)}, \quad f_2(\rho,q) = \frac{B_2(q) q \rho}{\left(1+\rho^2\right)
  \left(1+q^2\rho^2\right)},     \label{eq:f}
\end{equation}
in terms of two undetermined constants $B_1(q)$, and~$B_2(q)$.\par
%
%
Let us now further restrict the sesquilinear form~(\ref{eq:ansatz}) by imposing
that it is Hermitian:
 $ \overline{\langle \psi_1 | \psi_2 \rangle_q} = \langle \psi_2 | \psi_1
\rangle_q$.
We get the relation~\cite{nous}
\begin{equation}
 B_2(q) = \overline{B_1(q)}.
\label{eq:B}
\end{equation}
%
%
All the Hermiticity conditions on the su$_q$(2)
generators are satisfied by the form defined in Eqs.~(\ref{eq:ansatz}),
(\ref{eq:a}),
(\ref{eq:f}), (\ref{eq:B}), and the functions $\Psi^J_{MNq}(z,\zb)$, defined in
Eq.~(\ref{eq:Psi-z}), are  orthogonal with respect to such a form.
\par
%
%
To make $\langle \psi_1 | \psi_2 \rangle_q$ into a scalar product, it only
remains to
impose that it is a positive definite form,
which amounts to the condition
$ \left\langle \Psi^J_{JNq} | \Psi^J_{JNq} \right\rangle_q = 1.$
 A straighforward calculation of this squared norm leads to~\cite{nous}
\begin{equation}
  \frac{\ln q}{q-q^{-1}} \left(B_1(q)\, \overline{\gamma(J,N,q^{-1})} \,
\gamma(J,N,q)
  + \overline{B_1(q)}\, \overline{\gamma(J,N,q)} \,
\gamma(J,N,q^{-1})\right) = 1.
  \label{eq:norm-J-bis}
\end{equation}
We can choose
$\gamma(J,N,q) = 1$, and $ B_1(q) = (q-q^{-1})(2\ln q)^{-1}.$
When $q\to 1$, we find that $\langle
\psi_1 |
\psi_2 \rangle_q \to \langle \psi_1 | \psi_2 \rangle$, where the latter is given
by Eq.~(\ref{eq:prodscal}), as it should be.\par
%
%
 The results obtained can be summarized as follows:
%
\begin{proposition}   \label{prop-prodscal}
{}For $q \in \mbox{R}^+$, the scalar product
\begin{eqnarray}
  \langle \psi_1 | \psi_2 \rangle_q & = & \frac{q-q^{-1}}{2\ln q} \int dz\, d\zb
         \Biggl( \overline{\psi_1(z,\zb,q^{-1})}\, \frac{q^{-1}}{(1+z\zb)
(1+q^{-2}z\zb)}\,
         \psi_2(q^{-1}z,q^{-1}\zb,q) \nonumber \\
  & & \mbox{} + \overline{\psi_1(z,\zb,q)}\, \frac{q}{(1+z\zb) (1+q^2 z\zb)}\,
         \psi_2(qz,q\zb,q^{-1}) \Biggr),
\end{eqnarray}
unitarizes the su$_q$(2) realization~(\ref{eq:su-q}), where $N$ may take any
integer or half-integer value. The functions $\Psi^J_{MNq}(z,\zb)$, defined in
Eq.~(\ref{eq:Psi-z}), where $J = |N|$, $|N|+1$,~$\ldots$, $M = -J$, $-J+1$,
$\ldots$,~$J$, and $\gamma(J,N,q) = 1$, form an orthonormal set with respect to
such a scalar product.
\end{proposition}
%
\subsection{The case where \boldmath $q \in S^1$}
Whenever $q\in S^1$, the  scalar product obtained when $q \in \mbox{R}^+$
 does not work. In Ref.~\cite{nous}, by a treatment analogous to that of the
previous subsection, we established
\begin{proposition}  \label{prop-prodscalbis}
{}For generic $q \in S^1$, the scalar product
\begin{eqnarray}
  \langle \psi_1 | \psi_2 \rangle_q & = & \frac{q-q^{-1}}{2\ln q} \int dz\, d\zb
         \Biggl( \overline{\psi_1(z,\zb,q)}\, \frac{q^{-1}}{(1+z\zb)
(1+q^{-2}z\zb)}\,
         \psi_2(q^{-1}z,q^{-1}\zb,q) \nonumber \\
  & & \mbox{} + \overline{\psi_1(z,\zb,q^{-1})}\, \frac{q}{(1+z\zb) (1+q^2
z\zb)}\,
         \psi_2(qz,q\zb,q^{-1}) \Biggr)
\end{eqnarray}
unitarizes the su$_q$(2) realization~(\ref{eq:su-q}), where $N$ may take any
integer or half-integer value. The functions $\Psi^J_{MNq}(z,\zb)$, defined in
Eq.~(\ref{eq:Psi-z}), where $J = |N|$,
$|N|+1$,~$\ldots$, $M = -J$, $-J+1$, $\ldots$,~$J$, and $\gamma(J,N,q) =
1$, form
an orthonormal set with respect to such a scalar product.
\end{proposition}
%
%
\section{Conclusion}   \label{sec:conclusion}
In the present communication, we did extend the study of the su$_q$(2)
representations on the plane, carried out by Rideau and
Winternitz~\cite{rideau}, in two ways.

{}Firstly, we did prove that such
representations exist not only for $q \in
\mbox{R}^+$, but
also for generic $q \in S^1$. For such a purpose, we did provide an integral
representation for the functions $Q_{Jq}(\eta)$, entering the definition of the
$q$-Vilenkin functions, whenever $J$ takes any half-integer value.

Secondly, we did
unitarize the representations by determining appropriate scalar products
for both
ranges of $q$~values. Such scalar products are expressed in terms
of ordinary integrals, instead of $q$-integrals, as is usually the
case~\cite{vilenkin2}.

The resulting orthonormality relations for the $q$-Vilenkin
and related functions~\cite{nous}
should play an important role in applications to quantum mechanics, such as
those
considered in Ref.~\cite{irac1}.

%
%

%
\end{document}